\documentclass{amsart}
\usepackage{latexsym}
\usepackage{amssymb}
\usepackage{amsmath}

\begin{document}


\newtheorem{theorem}{Theorem}
\newtheorem{problem}{Problem}
\newtheorem{definition}{Definition}
\newtheorem{lemma}{Lemma}
\newtheorem{proposition}{Proposition}
\newtheorem{corollary}{Corollary}
\newtheorem{example}{Example}
\newtheorem{conjecture}{Conjecture}
\newtheorem{algorithm}{Algorithm}
\newtheorem{exercise}{Exercise}
\newtheorem{xample}{Example}
\newtheorem{remarkk}{Remark}

\newcommand{\be}{\begin{equation}}
\newcommand{\ee}{\end{equation}}
\newcommand{\bea}{\begin{eqnarray}}
\newcommand{\eea}{\end{eqnarray}}
\newcommand{\beq}[1]{\begin{equation}\label{#1}}
\newcommand{\eeq}{\end{equation}}
\newcommand{\beqn}[1]{\begin{eqnarray}\label{#1}}
\newcommand{\eeqn}{\end{eqnarray}}
\newcommand{\beaa}{\begin{eqnarray*}}
\newcommand{\eeaa}{\end{eqnarray*}}
\newcommand{\req}[1]{(\ref{#1})}

\newcommand{\lip}{\langle}
\newcommand{\rip}{\rangle}
\newcommand{\uu}{\underline}
\newcommand{\oo}{\overline}
\newcommand{\La}{\Lambda}
\newcommand{\la}{\lambda}
\newcommand{\eps}{\varepsilon}
\newcommand{\om}{\omega}
\newcommand{\Om}{\Omega}
\newcommand{\ga}{\gamma}
\newcommand{\ka}{\kappa}
\newcommand{\rrr}{{\Bigr)}}
\newcommand{\qqq}{{\Bigl\|}}

\newcommand{\dint}{\displaystyle\int}
\newcommand{\dsum}{\displaystyle\sum}
\newcommand{\dfr}{\displaystyle\frac}
\newcommand{\bige}{\mbox{\Large\it e}}
\newcommand{\integers}{{\Bbb Z}}
\newcommand{\rationals}{{\Bbb Q}}
\newcommand{\reals}{{\rm I\!R}}
\newcommand{\realsd}{\reals^d}
\newcommand{\realsn}{\reals^n}
\newcommand{\NN}{{\rm I\!N}}
\newcommand{\DD}{{\rm I\!D}}
\newcommand{\degree}{{\scriptscriptstyle \circ }}
\newcommand{\dfn}{\stackrel{\triangle}{=}}
\def\complex{\mathop{\raise .45ex\hbox{${\bf\scriptstyle{|}}$}
     \kern -0.40em {\rm \textstyle{C}}}\nolimits}
\def\hilbert{\mathop{\raise .21ex\hbox{$\bigcirc$}}\kern -1.005em {\rm\textstyle{H}}} 
\newcommand{\RAISE}{{\:\raisebox{.6ex}{$\scriptstyle{>}$}\raisebox{-.3ex}
           {$\scriptstyle{\!\!\!\!\!<}\:$}}} 

\newcommand{\hh}{{\:\raisebox{1.8ex}{$\scriptstyle{\degree}$}\raisebox{.0ex}
           {$\textstyle{\!\!\!\! H}$}}}

\newcommand{\OO}{\won}
\newcommand{\calA}{{\mathcal A}}
\newcommand{\calB}{{\mathcal B}}
\newcommand{\calC}{{\cal C}}
\newcommand{\calD}{{\cal D}}
\newcommand{\calE}{{\cal E}}
\newcommand{\calF}{{\mathcal F}}
\newcommand{\calG}{{\cal G}}
\newcommand{\calH}{{\cal H}}
\newcommand{\calK}{{\cal K}}
\newcommand{\calL}{{\mathcal L}}
\newcommand{\calM}{{\mathcal M}}
\newcommand{\calO}{{\cal O}}
\newcommand{\calP}{{\cal P}}
\newcommand{\calT}{{\mathcal T}} 
\newcommand{\calU}{{\mathcal U}}
\newcommand{\calX}{{\cal X}}
\newcommand{\calY}{{\mathcal Y}}
\newcommand{\calZ}{{\mathcal Z}}
\newcommand{\calXX}{{\cal X\mbox{\raisebox{.3ex}{$\!\!\!\!\!-$}}}}
\newcommand{\calXXX}{{\cal X\!\!\!\!\!-}}
\newcommand{\gi}{{\raisebox{.0ex}{$\scriptscriptstyle{\cal X}$}
\raisebox{.1ex} {$\scriptstyle{\!\!\!\!-}\:$}}}
\newcommand{\intsim}{\int_0^1\!\!\!\!\!\!\!\!\!\sim}
\newcommand{\intsimt}{\int_0^t\!\!\!\!\!\!\!\!\!\sim}
\newcommand{\pp}{{\partial}}
\newcommand{\al}{{\alpha}}
\newcommand{\sB}{{\cal B}}
\newcommand{\sL}{{\cal L}}
\newcommand{\sF}{{\cal F}}
\newcommand{\sE}{{\cal E}}
\newcommand{\sX}{{\cal X}}
\newcommand{\R}{{\rm I\!R}}
\renewcommand{\L}{{\rm I\!L}}
\newcommand{\vp}{\varphi}
\newcommand{\N}{{\rm I\!N}}
\def\ooo{\lip}
\def\ccc{\rip}
\newcommand{\ot}{\hat\otimes}
\newcommand{\rP}{{\Bbb P}}
\newcommand{\bfcdot}{{\mbox{\boldmath$\cdot$}}}

\renewcommand{\varrho}{{\ell}}
\newcommand{\dett}{{\textstyle{\det_2}}}
\newcommand{\sign}{{\mbox{\rm sign}}}
\newcommand{\TE}{{\rm TE}}
\newcommand{\TA}{{\rm TA}}
\newcommand{\E}{{\rm E\,}}
\newcommand{\won}{{\mbox{\bf 1}}}
\newcommand{\Lebn}{{\rm Leb}_n}
\newcommand{\Prob}{{\rm Prob\,}}
\newcommand{\sinc}{{\rm sinc\,}}
\newcommand{\ctg}{{\rm ctg\,}}
\newcommand{\loc}{{\rm loc}}
\newcommand{\trace}{{\,\,\rm trace\,\,}}
\newcommand{\Dom}{{\rm Dom}}
\newcommand{\ifff}{\mbox{\ if and only if\ }}
\newcommand{\nproof}{\noindent {\bf Proof:\ }}
\newcommand{\remark}{\noindent {\bf Remark:\ }}
\newcommand{\remarks}{\noindent {\bf Remarks:\ }}
\newcommand{\note}{\noindent {\bf Note:\ }}

\newcommand{\boldx}{{\bf x}}
\newcommand{\boldX}{{\bf X}}
\newcommand{\boldy}{{\bf y}}
\newcommand{\boldR}{{\bf R}}
\newcommand{\uux}{\uu{x}}
\newcommand{\uuY}{\uu{Y}}

\newcommand{\limn}{\lim_{n \rightarrow \infty}}
\newcommand{\limN}{\lim_{N \rightarrow \infty}}
\newcommand{\limr}{\lim_{r \rightarrow \infty}}
\newcommand{\limd}{\lim_{\delta \rightarrow \infty}}
\newcommand{\limM}{\lim_{M \rightarrow \infty}}
\newcommand{\limsupn}{\limsup_{n \rightarrow \infty}}

\newcommand{\ra}{ \rightarrow }

\newcommand{\ARROW}[1]
  {\begin{array}[t]{c}  \longrightarrow \\[-0.2cm] \textstyle{#1} \end{array} }

\newcommand{\AR}
 {\begin{array}[t]{c}
  \longrightarrow \\[-0.3cm]
  \scriptstyle {n\rightarrow \infty}
  \end{array}}

\newcommand{\pile}[2]
  {\left( \begin{array}{c}  {#1}\\[-0.2cm] {#2} \end{array} \right) }

\newcommand{\floor}[1]{\left\lfloor #1 \right\rfloor}

\newcommand{\mmbox}[1]{\mbox{\scriptsize{#1}}}

\newcommand{\ffrac}[2]
  {\left( \frac{#1}{#2} \right)}

\newcommand{\one}{\frac{1}{n}\:}
\newcommand{\half}{\frac{1}{2}\:}

\def\le{\leq}
\def\ge{\geq}
\def\lt{<}
\def\gt{>}

\def\squarebox#1{\hbox to #1{\hfill\vbox to #1{\vfill}}}
\newcommand{\nqed}{\hspace*{\fill}
          \vbox{\hrule\hbox{\vrule\squarebox{.667em}\vrule}\hrule}\bigskip}

\title{A sophisticated proof of the multiplication formula for
  multiple Wiener integrals}

\author{Ali  S\"uleyman  \"Ust\"unel}

\begin{abstract}
\noindent
We prove that  the formula which gives the Wiener chaos
decomposition of the  multiplication of two multiple Wiener integrals
with symmetric kernels is a straightforward application of the
Leibniz' formula.

\end{abstract}
\maketitle
\section{\bf{Introduction}}
\noindent
Let  $(W,H,\mu)$ be  the classical Wiener space, i.e.,
$W=C_0([0,1],\R^d)$, $H$ is the corresponding Cameron-Martin space
consisting of $\R^d$-valued  absolutely continuous functions
on $[0,1]$ with square integrable derivatives, which is a Hilbert
space under the norm $|h|_H^2=\int_0^1|\dot{h}(s)|^2ds$, where
$\dot{h}$ denotes the Radon-Nikodym derivative of the absolutely
continuous function $t\to h(t)$ w.r.t. the Lebesgue measure on $[0,1]$. Denote by  $(\calF_t,\,t\in
[0,1])$  the filtration of the canonical Wiener process, completed
w.r.t. $\mu$-negligeable sets. For $f\in H^{\hat{\otimes}p}$, i.e.,
the $p$-th order symmetric tensor product of $H$, we write $I_p(f)$
the multiple Wiener integral of $f$, which is defined as the iterated
integral:
$$
I_p(f)=p!\int_{t_1<t_2<\ldots t_p<1}f_p(t_1,\ldots,t_p) dW_{t_1}\ldots
dW_{t_p}\,.
$$
Let $g\in H^{\hat{\otimes}q}$, then the following multiplication
formula, which is due to Shigekawa (cf. \cite{Sh}) is extremely
important 
\begin{theorem}
\label{multiple}
We have
$$
I_p(f)I_q(g)=\sum_{i=0}^{p\wedge
  q}\frac{p!q!}{i!(p-i)!(q-i)!}I_{p+q-2i}(f\hat{\otimes}_ig)
$$
$\mu$- a.s., where $f\otimes_ig$ denotes the tensor $f\otimes g$ which
is contracted in its $2i$ components, i.e., 
$$
f\otimes_i g(t_1,\ldots,t_{p-i},s_1,\ldots,s_{q-i})=\int_{[0,1]^i}
f(t_1,\ldots,t_{p-i},u_1,\ldots,u_i)\,g(s_1,\ldots,s_{q-i},u_1,\ldots,u_i)
du_1\ldots du_i
$$
and $f\hat{\otimes}_ig$ is the symmetrization of $f\otimes_ig$ in its
remaining $p+q-2i$ variables.
\end{theorem}
\noindent
This theorem has been proved in 1980 by Shigekawa  using It\^o formula and
induction.  It is astonishing that no other proof has been seen in the
mathematical literature in spite of all the new techniques developped
thourough  the applications and the extensions of the Malliavin
calculus. 
We shall give here a completely new proof by relating
the formula given above to the Leibniz formula for $n$-th order
derivative of the multiplication of two smooth functions which is
given below,  whose proof
follows from its one-dimensional version:
\begin{lemma}
\label{Leibniz}
Assume that $F, \,G$ are two real-valued polynomials on $W$, then, for any $n\in \NN$,
we have 
$$
\nabla^n(F\,G)=\sum_{i=0}^n \binom{n}{i}\nabla^i
F\hat{\otimes}\nabla^{n-i}G\,.\,
$$
almost surely.
\end{lemma}
To make this note self-contained, we shall give also some results
connecting the Meyer distributions on Wiener space to the It\^o-Wiener
chaos decomposition of the elements of $L^2(\mu)$. For this we need
some notations which are explained in the next section.
\section{Notations}

We denote by $\nabla$ the Sobolev derivative on $(W,H,\mu)$ in the
direction of the Cameron-Martin space extended to $L^p(\mu)$, $p>1$,
the corresponding Sobolev spaces of real-valued dunctions are denoted
by $\DD_{p,k}$, $p>1,\,k\in \NN$, where $k$  denotes the degree of
differentiability and $p$ denotes the degree of invertibility.  For vector valued functions, we use
the notation $\DD_{p,k}(X)$, where $X$ is the range space. Note that,
for any $F\in \DD_{p,k}$, $\nabla F$ is an element of
$\DD_{p,k-1}(H)$. The formal adjoint of $\nabla$ w.r.to $\mu$ is
denoted as $\delta$ and called the divergence operator. It is easy to
see that $\DD_{p,1}(H)$ is in the domain of $\delta$ and for a $\xi\in
\DD_{p,1}(H)$, $\delta\xi$ coincides with the It\^o integral of
$\dot{\xi}$ if the latter is adapted to the Wiener filtration, here
$\dot{\xi}$ is the Sobolev derivative of $t\to\xi(t,w)\in H$. We
define the Ornstein-Uhlenbeck operator as $L=\delta\circ \nabla$, it
follows from Meyer inequalities that the seminorms defined by
$\|(I+L)^{k/2}F\|_{L^p(\mu,X)}$ are equivalent to the Sobolev norms
explained above for each $p>1$ and $k\in \NN$. Since the seminorms
defined with $L$ are also extendable to the case $k\in \R$ and $p>1$,
we obtain a scale of Banach spaces, still denoted by the same notation
$\DD_{p,k}(X)$, for $p>1$ and $k\in \R$. This construction implies
that $\nabla$ has a continuous extension as a map from $\DD'(X)\to
\DD'(X\otimes H)$, where 
$$
\DD'(X)=\cup_{p>1,k\in \R}\DD_{p.k}(X)
$$
and that $\delta$ has a continuous extension from $\DD'(X\otimes H)\to
\DD'(X)$, where the unions are all equipped with their inductive limit
topologies. We denote by $\DD(X)$ the intersection of the Sobolev
spaces $(\DD_{p,k}(X):\,p>1,k\in \R)$ equipped with the projective topology.

As an example of these considerations let $\phi\in \DD$ be a polynomial (i.e.,
$X=\R$), then for any $h\in H$, denoting by $\rho(\delta h)$ the Wick
exponential $\exp(\delta h-1/2|h|_H^2)$, due to Cameron-Martin
theorem,  we have 
\beaa
E[\phi\rho(\delta h)]&=&E[\phi(\cdot+h)]\\
&=&\sum_{n=0}^\infty
\frac{1}{n!}E[(\nabla^n\phi,h^{\otimes n})_{H^{\otimes n}}]\\
&=&\sum_{n=0}^\infty
\frac{1}{n!}(E[\nabla^n\phi],h^{\otimes n})_{H^{\otimes n}}]\\
&=&\sum_{n=0}^\infty
\frac{1}{(n!)^2}E\left[I_n(E[\nabla^n\phi]) I_n(h^{\otimes n})\right]\\
&=&\sum_{n=0}^\infty
\frac{1}{n!}E\left[I_n(E[\nabla^n\phi]) \rho(\delta h)\right]\,.
\eeaa
Since the linear combinations of the Wick exponentials are dense in
any $L^p(\mu),\,p>1$, we deduce from the above calculations, on the one
hand that 
\begin{equation}
\label{chaos-dec}
\phi=E[\phi]+\sum_{n=1}^\infty \frac{I_n(E[\nabla^n\phi])}{n!}
\end{equation}
$\mu$-a.s., and on the other hand that 
$$
\delta^n h^{\otimes n}=I_n(h^{\otimes n})
$$
$\mu$-a.s. By density of the linear combinations of $\{h^{\otimes
    n},\,h\in H\}$ in $H^{\hat{\otimes}n}$,  we deduce that 
$$
\delta^n \eta= I_n(\eta)
$$
$\mu$-a.s. for any $\eta\in H^{\hat{\otimes} n}$, where $\delta^n$
denotes the adjoint of $\nabla^n$ w.r.to the measure $\mu$. Furthermore
the  relation (\ref{chaos-dec}) has been given in an informal manner by Mc. Kean,
\cite{McK}, later D. Stroock has remarked that this expression extends to
$L^2$functionals since 
\beaa
(E[\nabla^n\phi],h^{\otimes n})_{H^{\otimes n}}&=&\langle
\nabla^n\phi,h^{\otimes n}\rangle\\
&=&\langle \phi,\delta^n h^{\otimes n}\rangle\,,
\eeaa
consequently
$$
|E[\nabla^n\phi]|_{H^{\otimes n}}\leq \sqrt{n!}\|\phi\|_{L^2(\mu)}\,.
$$
Let summarize what we have explained:
\begin{theorem}[McKean-Stroock]
The map defined on the smooth functions with values in $H^{\otimes
  n}$, defined as $\phi\to E[\nabla^n\phi]$ has a unique bounded
(linear) extension to the whole space $L^2(\mu)$ for any $n\geq 1$ and if we denote it
again with the same notation, then the following identity holds true
$$
\phi=E[\phi]+\sum_{i=1}^\infty \frac{1}{i!}\delta^i(E[\nabla^i\phi])
$$
where the sum converges in $L^2(\mu)$.
\end{theorem}
re
\section{Proof of the Multiplication formula}
Suppose that $p>q$ and let $\phi\in \DD$, using the identity
$\delta^pf=I_p(f)$ and the fact that $\delta^p$ is the adjoint of the
operator $\nabla^p$, we get, from Lemma \ref{Leibniz}
\beaa
E[I_p(f)I_q(g)\phi]&=&E[(f,\nabla^p(I_q(g)\phi)) _{H^{\hat{\otimes}p}}]\\
&=&E\left[\sum_{i=0}^p \binom{p}{i}(f,\nabla^i
I_q(g)\hat{\otimes}\nabla^{p-i}\phi) _{H^{\hat{\otimes}p}}\right]\\
&=&E\left[\sum_{i=0}^p \binom{p}{i} \frac{q!}{(q-i)!}                
 (f,I_{q-i}(g)\otimes\nabla^{p-i}\phi) _{H^{\otimes p}}\right]\\
&=&\sum_{i=0}^p
\binom{p}{i}\frac{q!}{(q-i)!} \,E\left[(f,\,I_{q-i}(g)\otimes\nabla^{p-i}\phi)_{H^{\otimes p}}\right]\\
&=&\sum_{i=0}^p
\binom{p}{i}\frac{q!}{(q-i)!} \,E\left[(I_{q-i}(g)\otimes_i\,f,\nabla^{p-i}\phi)_{H^{\otimes(p-i)}}\right]\\
&=&\sum_{i=0}^p
\binom{p}{i}\frac{q!}{(q-i)!} \,E[(g\otimes_i\,f,\nabla^{q-i}\nabla^{p-i}\phi)_{H^{\otimes(p+q-2i)}}]\\
&=&\sum_{i=0}^p
\binom{p}{i}\frac{q!}{(q-i)!} \,E[(g\otimes_i\,f,\nabla^{p+q-2i}\phi)_{H^{\otimes(p+q-2i)}}]\\
&=&\sum_{i=0}^p
\binom{p}{i}\frac{q!}{(q-i)!} \,E[I_{p+q-2i}(g\hat{\otimes}_i\,f),\phi]\,,
\eeaa
in the third equality we have used the fact that $f$ is a symmetric
tensor and the proof of Theorem \ref{multiple}  follows.
\qed


\vspace{2cm}

{\footnotesize{\bf{
\noindent
A.S. \"Ust\"unel, Telecom-Paristech (formerly ENST),
 Dept. Infres,\\
46, rue Barrault, 75013 Paris, France\\
email: ustunel@telecom-paristech.fr}
}}

\end{document}